# Large Deviations for the Fleming-Viot process with neutral mutation and selection[*]

D.A.Dawson  S. Feng

The Fields Institute  McMaster University

October 20, 1998

**Abstract**

Large deviation principles are established for the Fleming-Viot processes with neutral mutation and selection, and the corresponding equilibrium measures as the sampling rate goes to 0. All results are first proved for the finite allele model, and then generalized, through the projective limit technique, to the infinite allele model. Explicit expressions are obtained for the rate functions.

## 1 Introduction

The Fleming-Viot process is a measure-valued process describing the evolution of the distribution of genotypes in a population. In the case of two alleles it reduces to a one-dimensional diffusion process that approximates the classical Wright-Fisher model. The standard model in population genetics involves mutation, replacement sampling, and selective advantages among various genotypes.

Let $E$ be a compact metric space, $C(E)$ be the set of continuous functions on $E$, and $M_1(E)$ denote the space of all probability measures on $E$ with the topology of weak convergence. Let $A$ be the generator of a Markov process on $E$ with domain $D(A)$. Define $\mathcal{D} = \{F : F(\mu) = f(\langle \mu, \phi \rangle), f \in C_b^\infty(R), \phi \in D(A), \mu \in M_1(E)\}$. Then the generator of the Fleming-Viot process is

$$\mathcal{L}F(\mu) = \int_E \Big(A \frac{\delta F(\mu)}{\delta \mu(x)}\Big) \mu(dx) + \frac{\gamma}{2} \int_E \int_E \Big(\frac{\delta^2 F(\mu)}{\delta \mu(x) \delta \mu(y)}\Big) Q(\mu; dx, dy) \qquad (1.1)$$
$$= f'(\langle \mu, \phi \rangle) \langle \mu, A\phi \rangle + \frac{\gamma}{2} \int \int \phi(x) \phi(y) Q(\mu; dx, dy),$$

[*]Research supported by the Natural Science and Engineering Research Council of Canada

where

$$\delta F(\mu)/\delta\mu(x) = \lim_{\varepsilon \to 0+} \varepsilon^{-1}\{F(\mu + \varepsilon\delta_x) - F(\mu)\},$$
$$\delta^2 F(\mu)/\delta\mu(x)\delta\mu(y) = \lim_{\varepsilon_1 \to 0+,\, \varepsilon_2 \to 0+} (\varepsilon_1\varepsilon_2)^{-1}\{F(\mu + \varepsilon_1\delta_x + \varepsilon_2\delta_y) - F(\mu)\},$$
$$Q(\mu; dx, dy) = \mu(dx)\delta_x(dy) - \mu(dx)\mu(dy),$$

and $\delta_x$ stands for the Dirac measure at $x \in E$. The domain of $\mathcal{L}$ is $\mathcal{D}$. $E$ is called the type space, $A$ is known as the mutation operator, and the last term in (1.1) describes the continuous sampling. If the mutation operator has the form of $Af(x) = \frac{\theta}{2}\int(f(y) - f(x))\nu_0(dy)$ with $\nu_0 \in M_1(E)$, we call the process a Fleming-Viot process with neutral mutation. It is known that the Fleming-Viot process with neutral mutation has a unique reversible probability measure (cf. Ethier and Kurtz [6]).

In the present article we will consider the limiting behavior of this process as $\gamma \to 0$. In the first principal result we establish a large deviation principle (henceforth, LDP) for the sequence of reversible measures. It turns out that the sequence converges to the probability fixed point $\nu_0$ of the mutation operator exponentially fast and for $\theta = 1, \nu \ll \nu_0$ the deviation is characterized by the relative entropy defined as

$$H(\nu_0|\nu) = \begin{cases} \int_E h \log h\, d\nu & \text{if } \nu_0 \ll \nu \\ \infty & \text{otherwise}, \end{cases} \quad (1.2)$$

where $h$ is the Radon-Nikodym derivative of $\nu_0$ with respect to $\nu$. It is known, by Sanov's theorem, that the empirical measure $\{\frac{1}{n}\sum_{k=1}^{n}\delta_{X_k}\}_{n\geq 1}$ of an i.i.d. sequence of random variables with common distribution $\nu_0$ converges exponentially fast to $\nu_0$ as $n$ goes to infinity, and the deviation is characterized by the relative entropy $H(\nu|\nu_0)$. Our example here may be the first among the large deviation literature that has this "reversed" form of relative entropy as rate function. In Sanov's theorem, the influence of sampling is dominant while in the Fleming-Viot case this influence decreases to zero. This may be an explanation for the "reversed" expression of the two rate functions.

The second principal result of this article establishes a path level LDP for the Fleming-Viot process with neutral mutation and with selection. This can be viewed as an Freidlin-Wentzell type result in infinite dimension. Furthermore, the existing results on large deviations for finite dimensional diffusions usually assume either the diffusion coefficient is non-degenerate or the square root of the diffusion coefficient is uniformly Lipschitz continuous, but our model with finite alleles satisfies neither of them. Hence our results also include an extension of the finite dimensional Freidlin-Wentzell theory.

The sampling rate $\gamma$ can be interpreted as the inverse population size and the large deviation results of this paper describe the deviations from the "infinite population" deterministic limit.

The large deviation result for equilibrium measures is proved in Section 2 and the path level LDP is proved in Section 3. We will first prove the LDP for the Fleming-Viot process with neutral mutation, and then, by the Cameron-Martin-Girsanov transformation, for the case with selection. For processes without selection we will first prove the result for the finite allele model, and then generalize, through the projective limit technique, to the infinite allele model.

## 2 LDP for Equilibrium Measures

Let the following objects be given: $X$ a Hausdorff topological space, $\mathcal{B}_X$ a $\sigma$-algebra of space $X$, $\{P_\varepsilon\}_{\varepsilon>0}$ a family of probability measures on $(X, \mathcal{B}_X)$, $\{a_\varepsilon\}_{\varepsilon>0}$ a family of positive numbers tending to zero as $\varepsilon$ goes to zero, and a function $I : X \longrightarrow [0, \infty]$.

**Definition 2.1** *The function $I : X \longrightarrow [0, \infty]$ is called a rate function if it is lower semicontinuous. A rate function is called the good rate function if for any $r \geq 0$, the level set $\Phi_I(r) = \{x \in X : I(x) \leq r\}$ is compact. The constant $a_\varepsilon$ is called the speed.*

**Definition 2.2** *$\{P_\varepsilon\}$ satisfies a LDP with the rate function (or good rate function) $I$ if*

1. *for each $\mathcal{B}_X$-measurable open subset $G$ of $X$*

$$\liminf_{\varepsilon \to 0} a_\varepsilon \log P_\varepsilon(G) \geq - \inf_{x \in G} I(x); \tag{2.1}$$

2. *for each $\mathcal{B}_X$-measurable closed subset $B$ of $X$*

$$\limsup_{\varepsilon \to 0} a_\varepsilon \log P_\varepsilon(B) \leq - \inf_{x \in B} I(x); \tag{2.2}$$

**Remark:** The $\mathcal{B}_X$-measurable condition is needed when $\mathcal{B}$ is not the Borel $\sigma$-algebra and not all open or closed sets are $\mathcal{B}_X$-measurable. This situation may occur when the space $X$ is not separable. When the space $X$ is compact, all rate functions are good rate functions. Also the function $I$ in the above definition is unique when $X$ is a regular topological space. For an excellent introduction to basic concepts and techniques of large deviations refer to [4].

In this section we will establish LDP for the equilibrium measures of the Fleming-Viot processes with neutral mutation and with selection. We will start with the case when the type space $E$ is finite. Then, by using the projective limit technique, we obtain results for the case of $E = [0, 1]$.

For any $n \geq 1$, let $E = \{1, 2, \cdots, n\}$. The space $M_1(E)$, the set of all probability measures on $E$, can be identified with the $(n-1)$-dimensional simplex

$$\Delta_n = \{x = (x_1, \cdots, x_n) : x_i \geq 0, \ i = 1, \cdots, n; \ \sum_{i=1}^{n} x_i = 1\}.$$

The Fleming-Viot process with neutral mutation reduces to the neutral one-locus $n$-allele diffusion process with generator

$$\mathcal{A}_\gamma = \frac{\gamma}{2} \sum_{i,j=1}^{n} x_i(\delta_{ij} - x_j) \frac{\partial^2}{\partial x_i \partial x_j} + \sum_{i=1}^{n} \Big(\sum_{j=1}^{n} x_j q_{ji}\Big) \frac{\partial}{\partial x_i},$$

where $q_{ji}(j \neq i)$ is the intensity of a mutation from allele $j$ to allele $i$, and $q_{jj} = -\sum_{i \neq j} q_{ji}$.

Let $m$ denote the Lebesgue measure on $\Delta_n$. If the infinitesimal matrix $(q_{ij})$ is irreducible, then this diffusion has a unique stationary distribution which is absolutely continuous with respect to $m$. (cf. Shiga [11])

In the special case of parent-independent mutation, i.e.,

$$q_{ij} = \frac{\theta}{2} p_j > 0, \ 1 \leq i \neq j \leq n, \ \sum_{i=1}^{n} p_i = 1, \ \theta > 0.$$

Wright [14] discovered that for $p = (p_1, \cdots, p_n)$, the unique stationary distribution $\Pi_{\theta, \gamma, p} \in M_1(\Delta_n)$ is the Dirichlet distribution with parameters $p_1, \cdots, p_n$ given by

$$\Pi_{\theta, \gamma, p}(dx) = \frac{\Gamma[\gamma^{-1}\theta]}{\Gamma(\gamma^{-1}\theta p_1) \cdots \Gamma(\gamma^{-1}\theta p_n)} x_1^{\gamma^{-1}\theta p_1 - 1} \cdots x_n^{\gamma^{-1}\theta p_n - 1} dx_1 \cdots dx_{n-1}.$$

For any Borel measurable subset $C$ of $\Delta_n$, we have by Stirling's formula

$$\log \Pi_{\theta, \gamma, p}(C) = \log \Big\{ \frac{\sqrt{2\pi}(\gamma^{-1}\theta)^{\gamma^{-1}\theta - \frac{1}{2}} e^{\frac{\alpha \gamma}{12\theta}}}{(\sqrt{2\pi})^n (\gamma^{-1}\theta p_1)^{\gamma^{-1}\theta p_1 - \frac{1}{2}} e^{\frac{\alpha_1 \gamma}{12\theta p_1}} \cdots (\gamma^{-1}\theta p_n)^{\gamma^{-1}\theta p_n - \frac{1}{2}} e^{\frac{\alpha_n \gamma}{12\theta p_n}}} \quad (2.3)$$

$$\times \int_C x_1^{\gamma^{-1}\theta p_1 - 1} \cdots x_n^{\gamma^{-1}\theta p_n - 1} dx_1 \cdots dx_{n-1} \Big\}$$

$$= \frac{n-1}{2} \log \frac{1}{2\gamma\pi} + \frac{1}{2} \log \frac{(\theta p_1) \cdots (\theta p_n)}{\theta} + \frac{\gamma}{12\theta}\Big[\alpha - \frac{\alpha_1}{p_1} - \cdots - \frac{\alpha_n}{p_n}\Big]$$

$$- \gamma^{-1}\theta \sum_{i=1}^{n} p_i \log p_i + \log \int_C x_1^{\gamma^{-1}\theta p_1 - 1} \cdots x_n^{\gamma^{-1}\theta p_n - 1} dx_1 \cdots dx_{n-1},$$

where $0 < \alpha, \alpha_1, \cdots, \alpha_n < 1$ are some constants.

For any $\varepsilon > 0$, let $C_\varepsilon = \{x \in C : \min_{1 \leq i \leq n} x_i \geq \varepsilon\}$. For any measurable function $f$ on $\Delta_n$, $\|f\|_{L^{\gamma^{-1}}}$ denotes the $L^{\gamma^{-1}}$ norm of $f$ with respect to measure $m$. Then we have for $\gamma < \min_{1 \leq i \leq n}\{\theta p_i\}$,

$$\|I_C e^{\sum_{i=1}^{n} \theta p_i \log x_i}\|_{L^{\gamma^{-1}}} = \Big[\int_C x_1^{\gamma^{-1}\theta p_1} \cdots x_n^{\gamma^{-1}\theta p_n} dx_1 \cdots dx_{n-1}\Big]^\gamma \quad (2.4)$$

$$\leq [\int_C x_1^{\gamma^{-1}\theta p_1-1}\cdots x_n^{\gamma^{-1}\theta p_n-1}dx_1\cdots dx_{n-1}]^\gamma = ||I_C e^{\sum_{i=1}^n (\theta p_i-\gamma)\log x_i}||_{L^{\gamma^{-1}}}$$

$$\leq ||I_{A_\varepsilon} e^{\sum_{i=1}^n \theta p_i \log x_i}||_{L^{\gamma^{-1}}} (\frac{1}{\varepsilon})^{n\gamma} + m(C\setminus C_\varepsilon)$$

$$\leq ||I_C e^{\sum_{i=1}^n \theta p_i \log x_i}||_{L^{\gamma^{-1}}} (\frac{1}{\varepsilon})^{\gamma} + m(C\setminus C_\varepsilon).$$

Letting $\gamma \to 0$, then $\varepsilon \to 0$, we get

$$\lim_{\gamma\to 0}[\int_C x_1^{\gamma^{-1}\theta p_1-1}\cdots x_n^{\gamma^{-1}\theta p_n-1}dx_1\cdots dx_{n-1}]^\gamma = \operatorname{ess\,sup}\{I_C e^{\sum_{i=1}^n \theta p_i \log x_i} : x \in \Delta_n\}. \tag{2.5}$$

For any subset $B$ of $\Delta_n$, we have

$$\operatorname{ess\,sup}\{I_B e^{\sum_{i=1}^n \theta p_i \log x_i} : x \in \Delta_n\} \leq e^{-\inf_{x\in B} \sum_{i=1}^n \theta p_i \log \frac{1}{x_i}}. \tag{2.6}$$

On the other hand for any open subset $G$ of $\Delta_n$,

$$\operatorname{ess\,sup}\{I_G e^{\sum_{i=1}^n \theta p_i \log x_i} : x \in \Delta_n\} = \operatorname{ess\,sup}\{e^{\sum_{i=1}^n \theta p_i \log x_i} : x \in G\} \tag{2.7}$$
$$= e^{-\inf_{x\in G} \sum_{i=1}^n \theta p_i \log \frac{1}{x_i}}.$$

By (2.3),(2.5), and (2.6), we get that for any closed subset $B$ of $\Delta_n$

$$\limsup_{\gamma\to 0} \gamma \log \Pi_{\theta,\gamma,p}(B) \leq -\theta \inf_{x\in B}\{\sum_{i=1}^n (p_i \log \frac{1}{x_i} + p_i \log p_i)\}. \tag{2.8}$$

From (2.3), (2.5), and (2.7), we obtain that for any open subset $G$ of $\Delta_n$

$$\liminf_{\gamma\to 0} \gamma \log \Pi_{\theta,\gamma,p}(G) \geq -\theta \inf_{x\in G}\{\sum_{i=1}^n (p_i \log \frac{1}{x_i} + p_i \log p_i)\}. \tag{2.9}$$

Note that in the present situation the relative entropy $H(p|x)$ of $p$ with respect to $x$ is given by

$$H(p|x) = \sum_{i=1}^n p_i \log \frac{p_i}{x_i},$$

and is non-negative, continuous. This combined with the compactness of $\Delta_n$ implies that all level sets of $\{x \in \Delta_n : H(p|x) \leq r\}$ are compact. Thus we have proved the following theorem.

**Theorem 2.1** *When the mutation is parent independent, the family $\{\Pi_{\theta,\gamma,p}\}_{\gamma>0}$ satisfies a LDP on space $\Delta_n$ with the good rate function $I_p^\theta(x) = \theta H(p|x)$ as $\gamma$ goes to zero.*

**Remark:** It is well known that the relative entropy $H(x|p)$ is the rate function describing the large deviations of the empirical measure of an i.i.d. sequence of random variables with common distribution $p$. In the case of $\theta = 1$ the rate function in Theorem 2.1 has an "reversed" expression $H(p|x)$.

In Theorem 2.1, the probability measure $p$ is positive for every type $i$ in $E$. By making $x_i = 0$ whenever $p_i = 0$, the Dirichlet distribution can be defined for any $p \in \Delta_n$ as follows. Without loss of generality, let $p = (p_1, \cdots, p_r, 0, \cdots, 0), r < n, \tilde{p} = (p_1, \cdots, p_r)$. We define

$$\Pi_{\theta,\gamma,p} = \Pi_{\theta,\gamma,\tilde{p}} \times \delta_0^{\otimes n-r}. \tag{2.10}$$

For any $p, x \in \Delta_n$, define

$$I_p^\theta(x) = \begin{cases} \theta H(p|x) & \text{if } x \ll p \\ \infty & \text{otherwise,} \end{cases} \tag{2.11}$$

Then a generalized version of Theorem 2.1 can be stated as follows.

**Theorem 2.2** *For any $p \in \Delta_n$, the family $\{\Pi_{\theta,\gamma,p}\}_{\gamma>0}$ satisfies a LDP on space $\Delta_n$ with the good rate function $I_p^\theta(x)$ as $\gamma$ goes to zero.*

Next we consider the case when the type space $E = [0,1]$. Let $M_1([0,1])$ be the collection of all probability measures on $[0,1]$ with the topology of weak convergence. Another topology on $M_1([0,1])$ needed in the sequel is called the $\tau$ topology which is the smallest topology such that for each bounded measurable function $f$ on $[0,1]$, the map $\mu \to \int_0^1 f(z)\mu(dz)$ is continuous. The $\tau$-topology on $M_1([0,1])$ is stronger than the weak topology. But when the type space $E$ is finite, the two topologies coincide with each other. We use $M_1^\tau([0,1])$ to denote the space of all probability measures on $[0,1]$ equipped with the $\tau$ topology.

Note that for each $x \in [0,1]$, let $V_x = \{\mu \in M_1([0,1]) : \mu(x) > 3/4\}$. Then each $V_x$ is an nonempty open set in the $\tau$ topology and for different $x, y \in [0,1]$, $V_x$ and $V_y$ are disjoint. Since there are uncountable number of such open sets, the space $M_1^\tau([0,1])$ is not separable. The $\sigma$-algebra $\mathcal{B}^\tau$ on $M_1^\tau([0,1])$ is defined to be the smallest $\sigma$-algebra such that for every bounded, measurable function $f$ on $[0,1]$, the map $\mu \to \int_0^1 f(x)d\mu(x)$ is measurable. It is known (cf. [4]) that the Borel $\sigma$-algebra of space $M_1([0,1])$ is the same as $\mathcal{B}^\tau$. We use $\mathcal{B}$ to denote this common $\sigma$-algebra throughout the remainder of this section.

Let

$$\mathcal{P} = \{\{B_1, \cdots, B_r\} : r \geq 1, B_1, \cdots, B_r \text{ is a partition of } [0,1] \text{ by Borel measurable sets}\} \tag{2.12}$$

be the collection of all finite partitions of $[0,1]$. For any $\mu \in M_1([0,1]), \jmath = \{B_1, \cdots, B_r\} \in \mathcal{P}$, define $\pi_\jmath(\mu) = (\mu(B_1), \cdots, \mu(B_r))$.

For any $\mu, \nu \in M_1([0,1])$, let $H(\mu|\nu)$ be the relative entropy of $\mu$ with respect to $\nu$ defined in (1.2). Then it is known (cf. [5]) that

$$H(\mu|\nu) = \sup_{g \in C([0,1])} \{\int_0^1 g\,d\mu - \log \int_0^1 e^g d\nu\} = \sup_{g \in B([0,1])} \{\int_0^1 g\,d\mu - \log \int_0^1 e^g d\nu\}, \tag{2.13}$$

where $C([0,1]), B([0,1])$ are the sets of continuous functions and bounded measurable functions on $[0,1]$, respectively.

The following expression of the relative entropy will be repeatedly used in the sequel.

**Lemma 2.3** *For any $\mu, \nu \in M_1([0,1])$,*

$$H(\mu|\nu) = \sup_{J \in \mathcal{P}} H(\pi_J(\mu)|\pi_J(\nu)). \tag{2.14}$$

**Proof:** If $\mu$ is not absolutely continuous with respect to $\nu$, then both side of (2.14) are infinity. Now we assume that $h = \frac{d\mu}{d\nu}$ exists. Then for any $J = (B_1, \cdots, B_r) \in \mathcal{P}$, choose a function $g \in B([0,1])$ such that

$$g(z) = \sum_{i: \mu(B_i) > 0} \log \frac{\mu(B_i)}{\nu(B_i)} I_{B_i}(z).$$

By direct calculation we have

$$H(\pi_J(\mu)|\pi_J(\nu)) = \int_0^1 g(z) d\mu - \log \int_0^1 e^{g(z)} d\nu.$$

By (2.13), we get $H(\pi_J(\mu)|\pi_J(\nu)) \leq H(\mu|\nu)$ which implies that $\sup_{J \in \mathcal{P}} H(\pi_J(\mu)|\pi_J(\nu)) \leq H(\mu|\nu)$. On the other hand, for any $n \geq 1$, let

$$h_n(x) = \sum_{k=1}^{n2^n} (k-1)/2^n I_{B_k}(x) + n I_{A_n}(x),$$

where $B_k = \{z : (k-1)/2^n \leq h(z) < k/2^n\}$, $A_n = \{z : h(z) \geq n\}$. The fact that $h$ is integrable with respect to $\nu$ implies that $h_n$ is an increasing sequence of nonnegative simple functions converging to $h$ almost surely with respect to $\nu$. Let $\ell = \{B_1, \cdots, B_{n2^n}, A_n\}$. Then we have

$$\sup_{J \in \mathcal{P}} H(\pi_J(\mu)|\pi_J(\nu)) \geq H(\pi_\ell(\mu)|\pi_\ell(\nu)) \geq \int_0^1 h_n \log h_n d\nu.$$

Let $n$ go to infinity, we get (2.14) by (1.2) and the monotone convergence theorem. $\square$

Now we are ready to prove large deviation results for equilibrium measures of some Fleming-Viot processes. A Fleming-Viot process with neutral mutation has a generator of the following form:

$$\mathcal{L}F(\mu) = \int_E \Big(A \frac{\delta F(\mu)}{\delta \mu(z)}\Big) \mu(dz) + \frac{\gamma}{2} \int_E \int_E \Big(\frac{\delta^2 F(\mu)}{\delta \mu(z) \delta \mu(y)}\Big) Q(\mu; dz, dy), \tag{2.15}$$

where $E = [0,1]$, $\theta > 0$, $\nu_0 \in M_1([0,1])$, and

$$Af(z) = \frac{\theta}{2} \int_0^1 (f(y) - f(z)) \nu_0(dy),$$

$$Q(\mu; dz, dy) = \mu(dz) \delta_z(dy) - \mu(dz) \mu(dy).$$

For any symmetric bounded measurable function $V(z,y) \in B([0,1] \times [0,1])$, let

$$V(\mu) = \int_0^1 \int_0^1 V(z,y)\mu(dz)\mu(dy),$$

and

$$\langle \frac{\delta V(\mu)}{\delta \mu}, \frac{\delta F}{\delta \mu} \rangle = \int_0^1 \int_0^1 \int_0^1 \frac{\delta F}{\delta \mu}(z)[V(z,y) - V(y,w)]\mu(dz)\mu(dy)\mu(dw).$$

Then the generator of a Fleming-Viot process with neutral mutation and selection takes the form:

$$\mathcal{L}_V F(\mu) = \mathcal{L}F(\mu) + \langle \frac{\delta V(\mu)}{\delta \mu}, \frac{\delta F}{\delta \mu} \rangle, \tag{2.16}$$

where $V$ is called the fitness function.

It is known (cf. theorem 3.1 and theorem 3.2 in [6]) that the martingale problem associated with generators $\mathcal{L}$ and $\mathcal{L}_V$ are well-posed. Shiga [12] shows that the Fleming-Viot process with generator $\mathcal{L}$ has a unique, reversible stationary distribution $\Pi_{\theta,\gamma,\nu_0} \in M_1(M_1([0,1]))$, which is the distribution of a $M_1([0,1])$-valued random variable $\nu$ characterized by the property that whenever $r \geq 2$ and $B_1, \cdots, B_r$ is a partition of $[0,1]$, $(\nu(B_1), \cdots, \nu(B_r))$ has the Dirichlet distribution with parameters $\frac{1}{\gamma}\theta\nu_0(B_1), \cdots, \frac{1}{\gamma}\theta\nu_0(B_r)$. (Note that zero parameters can be removed to create a Dirichlet distribution on a simplex with dimension less than $r$.) The Fleming-Viot process with neutral mutation and selection also has a unique, reversible stationary distribution given by

$$\Pi_{\theta,\gamma,\nu_0,V}(d\mu) = Z^{-1} \exp[\frac{V(\mu)}{\gamma}]\Pi_{\theta,\gamma,\nu_0}(d\mu), \tag{2.17}$$

where $Z$ is the normalizing constant.(cf. Ethier and Kurtz [7]).

**Theorem 2.4** *The family* $\{\Pi_{\theta,\gamma,\nu_0}\}$ *satisfies a LDP on space* $M_1^\tau([0,1])$ *with the good rate function*

$$I_{\nu_0}^\theta(\mu) = \begin{cases} \theta H(\nu_0|\mu) & if \mu \ll \nu_0 \\ \infty & otherwise, \end{cases} \tag{2.18}$$

**Proof:** First note that the set of all finite partitions $\mathcal{P}$, partially ordered by $\jmath \succ \imath$ iff $\jmath$ is finer than $\imath$, is a partially ordered right-filtering set. For every $\jmath = (B_1, \cdots, B_r) \in \mathcal{P}$, let

$$\mathcal{X}_\jmath = \{x = (x_{B_1}, \cdots, x_{B_r}) : x_{B_i} \geq 0, \ i = 1, \cdots, r; \ \sum_{i=1}^r x_{B_i} = 1\}.$$

For any $\imath = (C_1, \cdots, C_l), \jmath = (B_1, \cdots, B_r) \in \mathcal{P}, \jmath \succ \imath$, define

$$\pi_{\imath\jmath} : \mathcal{X}_\jmath \to \mathcal{X}_\imath, \ (x_{B_1}, \cdots, x_{B_r}) \to (\sum_{B_k \subset C_1} x_{B_k}, \cdots, \sum_{B_k \subset C_l} x_{B_k}).$$

Let $\mathcal{X}$ represent the projective limit of the family $\{\mathcal{X}_\jmath, \pi_{\imath\jmath}, \imath, \jmath \in J, \succ\}$, $\pi_\jmath$ be the projective mapping. Then we have $M_1([0,1]) \subset \mathcal{X}$ by identifying $\mu \in M_1([0,1])$ with the projective limit of $\{(\mu(B_1), \cdots, \mu(B_r)) : \jmath = (B_1, \cdots, B_r) \in \mathcal{P}\}$.

On the other hand, any element $\mu$ of $\mathcal{X}$ can be viewed as a finitely additive measure on $[0,1]$. For any $f \in C([0,1])$ and $\mu \in \mathcal{X}$, we define the following "abstract integral" of $f$ with respect to $\mu$:

$$\langle \mu, f \rangle = \lim_{n \to \infty} \sum_{k=0}^{2^n-1} f(2^{-n}k) x_{B_k},$$

where $B_0 = [0, 2^{-n})$, $B_k = (2^{-n}k, 2^{-n}(k+1)]$ for $k = 1, \cdots, 2^n - 1$. The existence of this limit is guaranteed since $\sum_{k=0}^{2^n-1} f(2^{-n}k) x_{B_k}$ is a bounded Cauchy sequence. From this definition and the fact that any decreasing sequence $f_n \in C([0,1])$ which converges to zero pointwisely converges to zero uniformly (Dini's theorem), we conclude that $\langle \mu, f \rangle$ is a linear form on space $C([0,1])$ satisfying:

(a) $\langle \mu, f \rangle \geq 0$ for $f \geq 0$;

(b) $\langle \mu, f \rangle = 1$ for $f \equiv 1$;

(c) $\langle \mu, f \rangle \to 0$ as $f \downarrow 0$ pointwise.

By the Daniell-Stone theorem (see e.g. page 197 of Bauer [1]), $\mu$ is a probability measure on the Borel $\sigma$-algebra of $[0,1]$. The projective topology obtained is just the $\tau$-topology. Hence we have $\mathcal{X} = M_1^\tau([0,1])$. By Theorem 2.2 and Theorem 3.3 of Dawson and Gärtner [3], we have that the family $\{\Pi_{\theta,\gamma,\nu_0}\}$ satisfies a LDP on space $\mathcal{X}$ with the good rate function

$$I(\mu) = \sup\{I^\theta_{\pi_\jmath(\nu_0)}(\pi_\jmath(\mu)) : \jmath \in \mathcal{P}\}. \tag{2.19}$$

From (2.3), we see that $I(\mu) = I^\theta_{\nu_0}(\mu)$. $\square$

Now let

$$C(\theta, \nu_0, V) = \sup_{\mu \in M_1(E)} \{V(\mu) - I^\theta_{\nu_0}(\mu)\}, \quad I^{\theta,V}_{\nu_0}(\mu) = C(\theta, \nu_0, V) - [V(\mu) - I^\theta_{\nu_0}(\mu)].$$

Then we have the following LDP.

**Theorem 2.5** *The family $\{\Pi_{\theta,\gamma,\nu_0,V}\}$ satisfies a LDP on space $M_1([0,1])$ with the good rate function $I^{\theta,V}_{\nu_0}(\mu)$.*

**Proof:** Since $C(\theta, \nu_0, V)$ is a constant, $V(\mu)$ is bounded continuous in the $\tau$-topology, any level set of function $I^{\theta,V}_{\nu_0}(\mu)$ is a $\tau$-closed subset of a level set of function $I^\theta_{\nu_0}(\mu)$. Hence $I^{\theta,V}_{\nu_0}(\mu)$ is a good rate function in the $\tau$-topology, and thus in the weak topology too.

By using Varadhan's Lemma, we have for any $F(\mu) \in C(M_1([0,1]))$,

$$\lim_{\gamma \to 0} \gamma \log Z = \lim_{\gamma \to 0} \gamma \log \int_{M_1([0,1])} e^{V(\mu)/\gamma} \Pi_{\theta,\gamma,\nu_0}(d\mu) = C(\theta, \nu_0, V),$$

$$\lim_{\gamma \to 0} \gamma \log \int_{M_1([0,1])} e^{F(\mu)/\gamma} \Pi_{\theta,\gamma,\nu_0,V}(d\mu) = C(\theta, \nu_0, F+V) - C(\theta, \nu_0, V).$$

Since $M_1([0,1])$ is compact and thus the family $\Pi_{\theta,\gamma,\nu_0,V}$ is exponential tight, by Bryc's inverse Varadhan Lemma (cf. section 4.4 of Dembo and Zeitouni [4]), we get that the family $\{\Pi_{\theta,\gamma,\nu_0,V}\}$ satisfies a LDP on space $M_1([0,1])$ with the good rate function $I_{\nu_0}^{\theta,V}(\mu)$. $\square$

## 3 Path Level LDP

In this section we will establish three LDPs at the path level: LDP for finite type (or allele) model, LDP for the Fleming-Viot process with neutral mutation, and LDP for the Fleming-Viot process with neutral mutation and selection. The novelty of our result for finite type model is that the diffusion coefficient of the corresponding diffusion process is degenerate and the square root of the diffusion coefficient is non-Lipschitz. This is an extension of the Freidlin-Wentzell theory.

### 3.1 LDP For Finite Allele Model

Define

$$S_n = \{x = (x_1, \cdots, x_{n-1}) : x_i \geq 0,\ i = 1, \cdots, n-1;\ \sum_{i=1}^{n-1} x_i \leq 1\}.$$

The Fleming-Viot process with finite allele and neutral mutation is a finite dimensional diffusion process described by the following system of stochastic differential equations

$$dx_k^\varepsilon(t) = b_k(x^\varepsilon(t))dt + \varepsilon \sum_{l=1}^{n-1} \sigma_{kl}(x^\varepsilon(t))dB_l(t),\ 1 \leq k \leq n-1, \tag{3.1}$$

where $x^\varepsilon(t) = (x_1^\varepsilon(t), \cdots, x_{n-1}^\varepsilon(t))$, $b_k(x^\varepsilon(t)) = \frac{\theta}{2}(p_k - x_k^\varepsilon(t))$, and $\sigma(x^\varepsilon(t)) = (\sigma_{kl}(x^\varepsilon(t)))_{1 \leq k,l \leq n-1}$ is given by

$$\sigma(x^\varepsilon(t))\sigma'(x^\varepsilon(t)) = D(x^\varepsilon(t)) = (x_k^\varepsilon(t)(\delta_{kl} - x_l^\varepsilon(t)))_{1 \leq k,l \leq n-1},$$

and $\varepsilon^2 = \gamma$, $p_k = \nu_0(k) > 0$, $B_l(t)$, $1 \leq l \leq n-1$ are independent Brownian motions.

For a fixed $T > 0$ and $x \in S_n$, let $C([0,T], S_n)$ be the space of all $S_n$-valued continuous functions on $[0,T]$ endowed with the uniform topology, and $P_x^\varepsilon$ denote the law of $x^\varepsilon(\cdot)$ starting at $x$.

**Lemma 3.1** *The family $\{P_x^\varepsilon\}_{\varepsilon > 0}$ is exponentially tight on $C([0,T], S_n)$ for all $x \in S_n$.*

**Proof:** Let $c_1 \geq \sup_{x \in \mathcal{S}_n} ||b(x)||$, $c_2 = \sup_{||\theta||=1, x \in \mathcal{S}_n} \langle \theta, D(x)\theta \rangle > 0$. For any $a > 0$, choose $m_0 \geq 1, \varepsilon_0 > 0$ such that for $m \geq m_0, \varepsilon \leq \varepsilon_0$, we have $2(n-1)m^2 T \geq 1$ and $a\varepsilon^{-2} > 1$. Define

$$K_a = \cap_{m \geq m_0} \{x(\cdot) \in C([0,T], \mathcal{S}_n); \sup_{t,s \in [0,T], |t-s| < \gamma_m} |x(t) - x(s)| \leq \rho_m\},$$

where $\gamma_m = \frac{1}{m^2}$, $\rho_m = c_1 \gamma_m + \{2c_2(n-1)a[\frac{1}{m} + \gamma_m \ln \frac{2(n-1)T}{\gamma_m}]\}^{1/2}$. Let $b(x) = (b_1(x), \cdots, b_{n-1}(x))$. Then we have

$$P_x^\varepsilon\{x(\cdot) \in K_a^c\} \leq \sum_{m=m_0}^\infty P\{\sup_{t,s \in [0,T], |t-s| < \gamma_m} |x^\varepsilon(t) - x^\varepsilon(s)| > \rho_m\} \quad (3.2)$$

$$\leq \sum_{m=m_0}^\infty P\{\sup_{t,s \in [0,T], |t-s| < \gamma_m} |x^\varepsilon(t) - x^\varepsilon(s) - \int_s^t b(x^\varepsilon(\tau))d\tau| > \rho_m - c_1 \gamma_m\}$$

$$\leq \sum_{m=m_0}^\infty \frac{2(n-1)T}{\gamma_m} \exp\Big(-\frac{(\rho_m - c_1 \gamma_m)^2}{2c_2(n-1)\gamma_m \varepsilon^2}\Big)$$

$$= \sum_{m=m_0}^\infty \frac{2(n-1)T}{\gamma_m} \exp\Big(-\frac{a}{\gamma_m \varepsilon^2}(1/m + \gamma_m \ln(2(n-1)T/\gamma_m))\Big)$$

$$\leq \sum_{m=m_0}^\infty \frac{2(n-1)T}{\gamma_m} \Big(\frac{\gamma_m}{2(n-1)T}\Big)^{a/\varepsilon^2} \exp(-am/\varepsilon^2) \leq \frac{\exp(-\frac{a}{\varepsilon^2})}{1 - \exp(-\frac{a}{\varepsilon^2})},$$

where the third inequality is obtained by using Theorem 4.2.1 in Stroock and Varadhan [13]. Letting $\varepsilon$ go to zero, we get

$$\limsup_{\varepsilon \to 0} \varepsilon^2 \log P_x^\varepsilon\{K_a^c\} \leq -a. \quad (3.3)$$

The lemma follows since $K_a$ is a compact set in $C([0,T], \mathcal{S}_n)$. □

Let $\partial \mathcal{S}_n$ and $\mathcal{S}_n^\circ$ denote the boundary and interior of $\mathcal{S}_n$, respectively. By direct calculation we get $det(A(x)) = x_1 \cdots x_{n-1}(1 - \sum_{i=1}^{n-1} x_i)$. Thus for any $x \in \mathcal{S}_n^\circ$, $D(x)$ is invertible and the inverse is given by

$$D^{-1}(x) = \begin{pmatrix} \frac{1 - \sum_{1 \leq i \leq n-1, i \neq 1} x_i}{x_1 x_n} & \frac{1}{x_n} & \cdots & \frac{1}{x_n} \\ \frac{1}{x_n} & \frac{1 - \sum_{1 \leq i \leq n-1, i \neq 2} x_i}{x_2 x_n} & \cdots & \frac{1}{x_n} \\ \cdots & \cdots & \cdots & \cdots \\ \frac{1}{x_n} & \frac{1}{x_n} & \cdots & \frac{1 - \sum_{1 \leq i \leq n-1, i \neq n-1} x_i}{x_{n-1} x_n} \end{pmatrix} \quad (3.4)$$

where $x_n = 1 - \sum_{i=1}^{n-1} x_i$.

For any $x$ in $\mathcal{S}_n$ and $\varphi$ in $C([0,T], \mathcal{S}_n)$, let

$$H_1^x = \{\varphi : \varphi(t) = x + \int_0^t g(s)ds, g \in L^2([0,T])\},$$

$$K_\varphi = \{g \in H_1^x : \varphi(t) = x + \int_0^t b(\varphi(s))ds + \int_0^t \sigma(\varphi(s))\dot{g}(s)ds\},$$

and define

$$I_x(\varphi) = \begin{cases} \frac{1}{2}\inf_{g \in K_\varphi} \int_0^T (\dot{g}(t))^2 dt, & \varphi \in H_1^x \\ \infty, & \varphi \notin H_1^x \end{cases} \quad (3.5)$$

By using the explicit expression of $D^{-1}(x)$, we get the following for functions $\varphi(t)$ whose pathes are completely contained inside $\mathcal{S}_n^\circ$.

$$[\dot{\varphi}(t) - b(\varphi(t))] D^{-1}(\varphi(t))[\dot{\varphi}(t) - b(\varphi(t))]' \quad (3.6)$$

$$= \sum_{i,j=1, i \neq j}^{n-1} [\dot{\varphi}_i(t) - b_i(\varphi(t))][\dot{\varphi}_j(t) - b_j(\varphi(t))]\Big(1 - \sum_{k=1}^{n-1} \varphi_k(t)\Big)^{-1}$$

$$+ \sum_{i=1}^{n-1} [\dot{\varphi}_i(t) - b_i(\varphi(t))]^2 \Big(\frac{1 - \sum_{k=1, k \neq i}^{n-1} \varphi_k(t)}{\varphi_i(t)(1 - \sum_{k=1}^{n-1} \varphi_k(t))}\Big)$$

$$= \sum_{i=1}^n \frac{(\dot{\varphi}_i(t) - b_i(\varphi(t)))^2}{\varphi_i(t)},$$

where $\varphi_n(t) = 1 - \sum_{i=1}^{n-1} \varphi_i(t)$.

**Lemma 3.2** *If $\varphi$ hits the boundary $\partial \mathcal{S}_n$, then $I_x(\varphi) = \infty$.*

**Proof:** We will prove this result by contradiction. Assume there is a $\varphi$ such that $I_x(\varphi) < \infty$ and $\varphi$ hits the boundary. Let $t_0 \in (0, T]$ be the first time that $\varphi$ hits the boundary. Without loss of generality we further assume that the hitting occurs on the first coordinate $\varphi_1$. Now choose $0 < t_1 < t_2 < t_0$ such that $\inf_{t \in [t_1, t_2]} \{b_1(\varphi_1(t))\} > 0$, and $\log(\varphi_1(t))$ is absolutely continuous on $[t_1, t_2]$. By direct calculation we get

$$-2 \int_{t_1}^{t_2} \frac{\dot{\varphi}_1(t) b_1(\varphi(t))}{\varphi_1(t)} dt \leq \int_{t_1}^{t_2} \frac{(\dot{\varphi}_1(t) - b_1(\varphi(t)))^2}{\varphi_1(t)} dt \leq I_x(\varphi) < \infty.$$

On the other hand, performing integration by parts twice, we get

$$-2 \int_{t_1}^{t_2} \frac{\dot{\varphi}_1(t) b_1(\varphi(t))}{\varphi_1(t)} dt = -2\{b_1(\varphi(t_2)) \log(\varphi_1(t_2)) - b_1(\varphi(t_1)) \log(\varphi_1(t_1))$$

$$+ \frac{\theta}{2} \int_{t_1}^{t_2} \log(\varphi_1(t)) \dot{\varphi}_1(t) dt\}$$

$$= -2\{b_1(\varphi(t_2)) \log(\varphi_1(t_2)) - b_1(\varphi(t_1)) \log(\varphi_1(t_1))$$

$$+ \frac{\theta}{2}[\varphi_1(t_2)(\log(\varphi_1(t_2)) - 1) - \varphi_1(t_1)(\log(\varphi_1(t_1)) - 1)]\}.$$

Letting $t_2$ goes to $t_0$, the above will go to infinity because the first term goes to infinity while all other terms are bounded. This certainly contradicts the assumption of finiteness of $I_x(\varphi)$. □

Hence if we define $H_0^x = H_1^x \cap C([0,T], \mathcal{S}_n^\circ)$, then

$$I_x(\varphi) = \begin{cases} \frac{1}{2}\int_0^T \sum_{i=1}^n \frac{(\dot{\varphi}_i(t)-b_i(\varphi(t)))^2}{\varphi_i(t)}dt, & \varphi \in H_0^x \\ \infty, & \varphi \notin H_0^x \end{cases} \quad (3.7)$$

**Theorem 3.3** *For any $x \in \mathcal{S}_n^\circ$, the family $\{P_x^\varepsilon\}_{\varepsilon>0}$ satisfies a LDP on space $C([0,T], \mathcal{S}_n)$ with the good rate function $I_x(\cdot)$ and speed $\varepsilon^2$.*

**Remark:** If we introduce a map $\Psi$ between spaces $\mathcal{S}_n$ and $\Delta_n$ such that

$$\Psi(x_1, \cdots, x_{n-1}) = (x_1, \cdots, x_{n-1}, 1 - \sum_{i=1}^{n-1} x_i),$$

and for simplicity, let $\{P_x^\varepsilon\}_{\varepsilon>0}$ denote its image probability on space $C([0,T], \Delta_n)$ under the map $\Psi$, then by contraction principle we have that for any $x \in \Delta_n^\circ$, the family $\{P_x^\varepsilon\}_{\varepsilon>0}$ satisfies a large deviation principle on space $C([0,T], \Delta_n)$ with the good rate function $I_x(\cdot)$ and speed $\varepsilon^2$.

**Proof:** By Corollary 3.4 in Pukhalskii [10] and Lemma 3.1, it suffices to show that for every $\varphi \in C_x([0,T], \mathcal{S}_n)$,

$$\lim_{\delta \to 0} \liminf_{\varepsilon \to 0} \varepsilon^2 \log P_x^\varepsilon\{\sup_{t \in [0,T]} |x(t) - \varphi(t)| \leq \delta\} \quad (3.8)$$
$$= \lim_{\delta \to 0} \limsup_{\varepsilon \to 0} \varepsilon^2 \log P_x^\varepsilon\{\sup_{t \in [0,T]} |x(t) - \varphi(t)| \leq \delta\} = -I_x(\varphi).$$

First we assume that the path of $\varphi$ is contained in $\mathcal{S}_n^\circ$. Choose $\delta_0$ small enough such that

$$B \equiv \{x(\cdot) \in C_x([0,T], \mathcal{S}_n); \sup_{t \in [0,T]} |x(t) - \varphi(t)| \leq \delta_0\} \subset C_x([0,T], \mathcal{S}_n^\circ).$$

Then we have

$$\alpha \equiv \frac{1}{2} \inf_{t \in [0,T], x(\cdot) \in B} d(x(t), \partial \mathcal{S}_n) > 0.$$

where $d(y, \partial \mathcal{S}_n) = \inf_{z \in \partial \mathcal{S}_n} |y - z|$.

Define

$$\tilde{\sigma}(x) \equiv \begin{cases} \sigma(x), & d(x, \partial \mathcal{S}_n) \geq \alpha \\ \text{smooth and uniform Lipschitz continuous}, & \text{else} \end{cases}$$

By replacing $\sigma(x)$ with $\tilde{\sigma}(x)$, we define the diffusion $\tilde{x}^\varepsilon(t)$, the law $\tilde{P}_x^\varepsilon$, the matrix $\tilde{A}(x)$, and the function $\tilde{I}_x(\cdot)$ respectively. It is easy to see that $I_x(\varphi) = \tilde{I}_x(\varphi)$. By using Theorem 5.6.7 in Dembo and Zeitouni [4] and Theorem 3.5 of Chapter 3 in Freidlin and Wentzell [8], we have for any $f \in C_x([0,T], \mathcal{S}_n)$,

$$\lim_{\delta \to 0} \liminf_{\varepsilon \to 0} \varepsilon^2 \log \tilde{P}_x^\varepsilon \{ \sup_{t \in [0,T]} |x(t) - f(t)| \leq \delta \} \qquad (3.9)$$
$$= \lim_{\delta \to 0} \limsup_{\varepsilon \to 0} \varepsilon^2 \log \tilde{P}_x^\varepsilon \{ \sup_{t \in [0,T]} |x(t) - f(t)| \leq \delta \} = -\tilde{I}_x(f).$$

Replacing $f$ with $\varphi$ in (3.9), we end up with (3.8) because the corresponding probabilities in both equations are the same for $\delta < \delta_0$.

Now we are ready to prove (3.8) for paths that hit the boundary of $\mathcal{S}_n$. Let $\varphi$ be such a path. The following is trivially true.

$$\lim_{\delta \to 0} \liminf_{\varepsilon \to 0} \varepsilon^2 \log P_x^\varepsilon \{ \sup_{t \in [0,T]} |x(t) - \varphi(t)| \leq \delta \} \geq -I_x(\varphi) = -\infty. \qquad (3.10)$$

On the other hand, let $t_0 > 0$ be the first time when $\varphi$ hits the boundary. Then for any $t \in (0, t_0)$, $\varphi$ will not hit the boundary on $[0, t]$. By using an argument similar to that used in the derivation of (3.9) we get

$$\lim_{\delta \to 0} \limsup_{\varepsilon \to 0} \varepsilon^2 \log \tilde{P}_x^\varepsilon \{ \sup_{s \in [0,t]} |x(s) - \varphi(s)| \leq \delta \} = -I_x^t(\varphi), \qquad (3.11)$$

where $I_x^t(\varphi)$ is the restriction of $I_x(\varphi)$ on $[0, t]$. Hence

$$\lim_{\delta \to 0} \limsup_{\varepsilon \to 0} \varepsilon^2 \log \tilde{P}_x^\varepsilon \{ \sup_{s \in [0,T]} |x(s) - \varphi(s)| \leq \delta \} \qquad (3.12)$$
$$\leq \lim_{\delta \to 0} \limsup_{\varepsilon \to 0} \varepsilon^2 \log \tilde{P}_x^\varepsilon \{ \sup_{s \in [0,t]} |x(s) - \varphi(s)| \leq \delta \} = -I_x^t(\varphi).$$

From the proof of lemma 3.2, we get that $I_x^t(\varphi)$ converges to infinity as $t \nearrow t_0$ which implies

$$\lim_{\delta \to 0} \limsup_{\varepsilon \to 0} \varepsilon^2 \log \tilde{P}_x^\varepsilon \{ \sup_{s \in [0,T]} |x(s) - \varphi(s)| \leq \delta \} = -\infty. \qquad (3.13)$$

Lemma 3.2 combined with (3.10) and (3.13) implies that (3.8) holds for all $\varphi \in C([0,T], \mathcal{S}_n)$. $\square$

## 3.2 LDP for Fleming-Viot Processes with Neutral Mutation

The Fleming-Viot process with neutral mutation has many nice properties. One of them is called the partition property, namely, given any finite partition of the type space $E = \cup_{i=1}^K E_i$, then $\{X(t, E_i) : i = 1, \cdots, K\}$ is a finite dimensional diffusion process as in section 3.1.(cf. Ethier and Kurtz [7].) By using this property and the projective limit technique, we will establish a LDP for the Fleming-Viot process with neutral mutation. This kind of result can be viewed as an infinite dimensional generalization of the Freidlin-Wentzell theory. In the remainder of this article we will have $E = [0, 1]$.

Consider the following family of partitions of space $E$

$$\mathcal{J} = \{\{[0, t_1], (t_1, t_2], \cdots, (t_n, 1]\} : 0 < t_1 < t_2 < \cdots < t_n < 1; n = 1, 2, \cdots\},$$

with the same partial ordering $\succ$ as in $\mathcal{P}$ of (2.12). It is clear that $\mathcal{J} \subset \mathcal{P}$. Thus we will denote generic element of $\mathcal{J}$ by $\imath, \jmath$, etc., and for any $\imath, \jmath \in \mathcal{J}$, the space $\mathcal{X}_\imath$ and the mappings $\pi_{\imath\jmath}, \pi_\imath$ are defined accordingly.

By an argument similar to that used in Section 2, we can show that the projective limit of $(\mathcal{X}_\imath, \pi_{\imath\jmath})_{\imath,\jmath \in \mathcal{J}}$ is $M_1(E)$ and projective topology is stronger than the weak topology but weaker than the $\tau$-topology. We will use $M_1^{pro}(E)$ to denote the space $M_1(E)$ with this projective limit topology.

Let $\mathcal{C}_\jmath = C([0, T], \mathcal{X}_\jmath)$ be equipped with the usual uniform topology. For any $\jmath = \{B_1, \cdots, B_r\} \succ \imath = \{C_1, \cdots, C_l\}$, define a map $p_{\imath\jmath}$ between spaces $\mathcal{C}_\jmath$ and $\mathcal{C}_\imath$ such that

$$p_{\imath\jmath} : \mathcal{C}_\jmath \to \mathcal{C}_\imath, \ (x_{B_1}(t), \cdots, x_{B_r}(t)) \to (\sum_{B_k \subset C_1} x_{B_k}(t), \cdots, \sum_{B_k \subset C_l} x_{B_k}(t)).$$

It is clear that $p_{\imath\jmath}$ is continuous, and for any $\ell \succ \jmath \succ \imath$, $p_{\imath\ell} = p_{\imath\jmath} \circ p_{\jmath\ell}$. Now let $\mathcal{C}$ be the projective limit of the family $\{(\mathcal{C}_\imath, p_{\imath\jmath}); \imath, \jmath \in \mathcal{J}\}$ and $p_\jmath : \mathcal{C} \to \mathcal{C}_\jmath$ be the corresponding projection. Obviously $C([0, T], M_1^{pro}([0, 1]))$ is a subset of $\mathcal{C}$. On the other hand, let $\{\mu^\jmath(\cdot) : \jmath \in \mathcal{J}\}$ be any element of $\mathcal{C}$. Then for any $t \in [0, T]$, $\{\mu^\jmath(t) : \jmath \in \mathcal{J}\}$ can be identified as a unique element $\mu(t)$ of $M_1^{pro}(E)$ and thus $\{\mu^\jmath(\cdot) : \jmath \in \mathcal{J}\}$ can be identified as $\mu(\cdot)$. For any $(a, b] \subset E$, by the definition of projective limit, $\mu(t)((a, b])$ is continuous in $t$. Hence $\mu(\cdot) \in C([0, T], M_1^{pro}(E))$ and $\mathcal{C}$ can be identified as $C([0, T], M_1^{pro}(E))$.

For any $\nu \in M_1(E)$, let $P_\nu^{\theta, \gamma, \nu_0}$ be the unique solution of the martingale problem associated with generator $\mathcal{L}$ in (2.15) starting at $\nu$.

**Definition 3.1** *A probability measure $\nu \in M_1(E)$ is called non-degenerate if for any $\jmath \in \mathcal{J}$, $\pi_\jmath(\nu)$ has no zero component. A path $\mu(\cdot) \in C([0, T], M_1^{pro}(E))$ is called non-degenerate if for every $t$ in $[0, T]$, $\mu(t)$ is non-degenerate.*

**Remark:** Any probability measure with support $E$ is non-degenerate.

**Definition 3.2** *A path $\mu(\cdot) \in C([0, T], M_1^{pro}(E))$ is called absolutely continuous if for every $\imath \in \mathcal{J}$, $\pi_\imath(\mu)(t)$ is absolutely continuous as a multidimensional real valued function.*

**Theorem 3.4** *For any non-degenerate $\nu \in M_1(E)$, the family $\{P_\nu^{\theta,\gamma,\nu_0}\}$ satisfies a LDP on space $C([0,T], M_1^{pro}(E))$ as $\gamma$ goes to zero with the good rate function*

$$I_\nu(\mu(\cdot)) = \sup_{\jmath \in \mathcal{J}} I_{\pi_\jmath(\nu)}(p_\jmath(\mu(\cdot))). \tag{3.14}$$

**Remark:** Let $\mu(t)$ be a path of the Fleming-Viot process with neutral mutation. Then for every $t > 0$, the support of $\mu(t)$ is a subset of the support of $\nu_0$. Therefore the essential part of the type space is the support of $\nu_0$. If we choose $E$ to be the support of $\nu_0$, then the result still holds. Because of this we assume in the sequel, without loss of generality, that the support of $\nu_0$ is $E = [0,1]$.

**Proof:** Since the path of Fleming-Viot process is continuous in the $\tau$-topology (cf. Shiga [12]), we have that $P_\nu^{\theta,\gamma,\nu_0}\{C([0,T], M_1^{pro}(E))\} = 1$. On the other hand for any $\jmath \in \mathcal{J}$, by the nondegeneracy of $\nu$, Theorem 3.3 and the remark following it, we have that $p_\jmath(P_\nu^{\theta,\gamma,\nu_0})$ satisfies a LDP on space $\mathcal{C}_\jmath$ with the good rate function $I_{\pi_\jmath(\nu)}(\cdot)$. Applying Theorem 3.3 in Dawson and Gärtner [3] we get the result. $\square$

Let $C^{1,0}([0,T] \times E)$ denote the set of all continuous functions on $[0,T] \times E$ with continuous first order derivative in time. For any $\nu \in M_1(E)$, $\mu(\cdot) \in C([0,T], M_1^{pro}(E))$, define

$$S_\nu(\mu(\cdot)) = \sup_{f \in C^{1,0}([0,T] \times E)} \{\langle \mu(T), f(T) \rangle - \langle \mu(0), f(0) \rangle \tag{3.15}$$
$$- \int_0^T \langle \mu(s), (\frac{\partial}{\partial s} + A)f \rangle \, ds - \frac{1}{2} \int_0^T \int \int f(s,x)\, f(s,y)\, Q(\mu(s); dx, dy)\, ds\},$$

and

$$\mathcal{H}_\nu = \{\mu(\cdot) \in C([0,T], M_1^{pro}(E)) : \mu(0) = \nu, \mu(\cdot) \text{ is absolutely continuous,}$$
$$\text{non-degenerate, and for any } \jmath \in \mathcal{J},\, p_\jmath(\mu)(\cdot) \in L^2([0,T])\}.$$

**Remark:** It is clear that $I_\nu(\mu(\cdot)) = \infty$ for any $\mu(\cdot)$ not in $\mathcal{H}_\nu$.

The next theorem gives an variational form of the rate function $I_\nu(\mu(\cdot))$.

**Theorem 3.5** *For any non-degenerate $\nu \in M_1(E)$, and any $\mu(\cdot) \in \mathcal{H}_\nu$, we have*

$$I_\nu(\mu(\cdot)) = S_\nu(\mu(\cdot)).$$

**Proof:** For any

$$\jmath = \{B_1 = [0, b_1), B_2 = [b_1, b_2), \cdots, B_{n-1} = [b_{n-1}, 1]\} \in \mathcal{J},$$

and $f(t,y) \in C^{1,0}([0,T] \times E)$, let $b_0 = 0$ and

$$\pi_\jmath(f)(t,y) = f(t, b_k), \text{ for } y \in B_k, k = 0, \cdots, n-1.$$

Then from the absolute continuity of $\mu(\cdot)$, we have

$$\int_0^T \langle \dot{\mu}(s), \pi_J(f)(s)\rangle\, ds = \langle \mu(T), \pi_J(f)(T)\rangle - \langle \mu(0), \pi_J(f)(0)\rangle - \int_0^T \langle \mu(s), \frac{\partial}{\partial s}\pi_J(f)(s)\rangle\, ds. \quad (3.16)$$

Approximating $f$ by $\pi_J(f)$, we obtain the following result of integration by parts:

$$\int_0^T \langle \dot{\mu}(s), f(s)\rangle\, ds = \langle \mu(T), f(T)\rangle - \langle \mu(0), f(0)\rangle - \int_0^T \langle \mu(s), \dot{f}(s)\rangle\, ds. \quad (3.17)$$

Let $x = \pi_J(\nu)$. For $k = 0, \cdots, n-1$, let $\varphi_k(t) = p_J(\mu)(t)(B_k)$, $\varphi(t) = (\varphi_0(t), \cdots, \varphi_{n-1}(t))$ and $g(t, y) = (\dot{\varphi}_k(t) - b_k(\varphi(t)))^2/\varphi_k(t)$ for $y \in B_k$. Here ⌈corrupted⌉$g(t,y)$ is well defined since $\mu(\cdot)$ is non-degenerate, and absolutely continuous. Then by direct calculation we have

$$\begin{aligned}
I_{\pi_J(\nu)}(p_J(\mu(\cdot))) = I_x(\varphi) &= \frac{1}{2}\sum_{k=0}^{n-1}\int_0^T \frac{(\dot{\varphi}_k(s) - b_k(\varphi(s)))^2}{\varphi_k(s)}\,ds \\
&= \int_0^T \langle \dot{\mu}(s), g(s)\rangle\,ds - \int_0^T \langle \mu(s), Ag\rangle\,ds \\
&\quad - \frac{1}{2}\int_0^T \int_E\int_E g(s,x)\,g(s,y)\,Q(\mu(s);dx,dy)\,ds \quad (3.18)\\
&\leq \sup_{f \in C^{1,0}([0,T]\times E)}\{\int_0^T \langle \dot{\mu}(s), \pi_J(f)(s)\rangle\,ds - \int_0^T \langle \mu(s), A(\pi_J(f))\rangle\,ds \\
&\quad - \frac{1}{2}\int_0^T \int_E\int_E \pi_J(f)(s,x)\,\pi_J(f)(s,y)\,Q(\mu(s);dx,dy)\,ds\},
\end{aligned}$$

where the inequality holds because $g$ can be approximated pointwise by functions in the set $\{\pi_J(f) : f \in C^{1,0}([0,T]\times E)\}$.

On the other hand,

$$\begin{aligned}
\sup_{f \in C^{1,0}([0,T]\times E)}&\{\int_0^T \langle \dot{\mu}(s), \pi_J(f)(s)\rangle\,ds - \int_0^T \langle \mu(s), A(\pi_J(f))\rangle\,ds \\
&\quad - \frac{1}{2}\int_0^T \int_E\int_E \pi_J(f)(s,x)\,\pi_J(f)(s,y)\,Q(\mu(s);dx,dy)\,ds\} \\
&\leq \sup_{f \in C^{1,0}([0,T]\times E)}\int_0^T [\sum_{k=1}^{n-1} f(s,B_k)(\dot{\varphi}_k(s) - b_k(\varphi(s))) \quad (3.19)\\
&\quad - \frac{1}{2}\sum_{k,l=1}^{n-1} f(s,B_k)D_{kl}(\varphi(s))f(s,B_l)]\,ds \\
&\leq \int_0^T \sup_{\alpha \in R^n}[\sum_{k=1}^{n-1}\alpha_k(\dot{\varphi}_k(s) - b_k(\varphi(s))) - \frac{1}{2}\sum_{k,l=1}^{n-1}\alpha_k D_{kl}(\varphi(s))\alpha_l]\,ds \\
&= I_x(\varphi),
\end{aligned}$$

which, together with (3.18), implies that

$$I_{\pi_J(\nu)}(p_J(\mu(\cdot))) = \sup_{f \in C^{1,0}([0,T]\times E)}\{\int_0^T \langle \dot{\mu}(s), \pi_J(f)(s)\rangle\,ds - \int_0^T \langle \mu(s), A(\pi_J(f))\rangle\,ds \quad (3.20)$$

$$-\frac{1}{2}\int_0^T\int_E\int_E \pi_{\jmath}(f)(s,x)\,\pi_{\jmath}(f)(s,y)\,Q(\mu(s);dx,dy)\,ds\}.$$

By using the inequality

$$\frac{(a+b)^2}{c+d} \leq \frac{a^2}{c} + \frac{b^2}{d}, a,b,c,d > 0,$$

and the expression (3.7) for $I_x(\varphi)$, we have that for any $\jmath_2 \succ \jmath_1$,

$$I_{\pi_{\jmath_2}(\nu)}(p_{\jmath_2}(\mu(\cdot))) \geq I_{\pi_{\jmath_1}(\nu)}(p_{\jmath_1}(\mu(\cdot))).$$

Taking the supremum on both sides of (3.20) over the set $\mathcal{J}$, we finally get

$$\begin{aligned} I_{(\nu)}((\mu(\cdot)) &= \sup_{f\in C^{1,0}([0,T]\times E)} \{\int_0^T \langle \dot{\mu}(s), f(s)\rangle ds - \int_0^T \langle \mu(s), Af\rangle\,ds \\ &\quad -\frac{1}{2}\int_0^T\int_E\int_E f(s,x)\,f(s,y)\,Q(\mu(s);dx,dy)\,ds\} \\ &= S_\nu(\mu(\cdot)), \end{aligned} \qquad (3.21)$$

where the last equality follows from integration by parts (3.17). $\square$

Let $C'^\infty(E)$ be the family of all continuous functions on $E$ possessing continuous derivatives of all order. For any linear functional $\vartheta$ on space $C'^\infty(E)$, define

$$||\vartheta||_\mu^2 = \sup_{f\in C^\infty(E)} [\langle \vartheta, f\rangle - \frac{1}{2}\int_E\int_E f(x)\,f(y)\,Q(\mu;dx,dy)].$$

Then we have

**Theorem 3.6** *For any non-degenerate $\nu \in M_1(E)$, and any $\mu(\cdot)$ in $\mathcal{H}_\nu$, we have*

$$I_\nu(\mu(\cdot)) = \int_0^T ||\dot{\mu}(s) - A^*(\mu(s))||_{\mu(s)}^2\,ds, \qquad (3.22)$$

*where $A^*$ is the formal adjoint of $A$ define through the equality $\langle A^*(\mu), f\rangle = \langle \mu, Af\rangle$.*

We defer the proof of this theorem to Appendix.

### 3.3 LDP for Fleming-Viot Processes with Selection

Finally we turn to prove the LDP of the Fleming-Viot process with selection. The generator of the process is given in (2.16). We will assume that the fitness function $V(x,y)$ is continuous on $E^{\otimes 2}$ in the sequel.

For any $\nu \in M_1(E)$, let $P_\nu^{\theta,\gamma,V,\nu_0}$ be the unique solution of the martingale problem associated with generator $\mathcal{L}_V$ started at $\nu$. For simplicity we will not distinguish between $P_\nu^{\theta,\gamma,\nu_0}$, $P_\nu^{\theta,\gamma,V,\nu_0}$ and their respective restrictions on $C([0,T], M_1^{pro}(E))$. By the Cameron-Martin-Girsanov transformation (see Dawson [2]) we have that, restricted on $C([0,T], M_1^{pro}(E))$,

$$\frac{dP_\nu^{\theta,\gamma,\nu_0}}{dP_\nu^{\theta,\gamma,V,\nu_0}} = Z_V(T) = \exp[\frac{1}{\gamma} G_V(\mu(\cdot))] > 0, \tag{3.23}$$

where

$$G_V(\mu(\cdot)) = \int_0^T \int_E [\int_E V(y,z)\mu(s,dz)] M(ds,dy) \tag{3.24}$$
$$- \frac{1}{2} \int_0^T \int_E \int_E [\int_E V(x,z)\mu(s,dz)][\int_E V(y,z)\mu(s,dz)] Q(\mu(s); dx, dy) ds,$$

and $M(ds, dy)$ is the martingale measure obtained from the martingale

$$M_t(\phi) = \langle \mu(t), \phi \rangle - \langle \mu(0), \phi \rangle - \int_0^t \langle \mu(s), A\phi \rangle ds.$$

Let $\rho$ be the Prohorov metric on $M_1(E)$, and $C^w([0,T], M_1^{pro}(E))$ denote $C([0,T], M_1^{pro}(E))$ equipped with the subspace topology of $C([0,T], M_1(E))$. For any $v(\cdot), \mu(\cdot) \in C^w([0,T], M_1^{pro}(E))$, let

$$\varrho(v(\cdot), \mu(\cdot)) = \sup_{t \in [0,T]} \rho(v(t), \mu(t)).$$

Define

$$R(\mu, dx) = \int_E [\int_E V(y,z)\mu(dz)] Q(\mu; dx, dy).$$

Then we have the following theorem.

**Theorem 3.7** *For any non-degenerate $\nu \in M_1(E)$, the family $\{P_\nu^{\theta,\gamma,V,\nu_0}\}$ satisfies the following local LDP on space $C^w([0,T], M_1^{pro}(E))$ as $\gamma$ goes to zero with the good rate function $I_{\nu,V}(\mu(\cdot))$:*

$$\lim_{\delta \to 0} \liminf_{\gamma \to 0} \gamma \log P_\nu^{\theta,\gamma,V,\nu_0}\{\varrho(v(\cdot), \mu(\cdot)) < \delta\} \tag{3.25}$$
$$= \lim_{\delta \to 0} \limsup_{\gamma \to 0} \gamma \log P_\nu^{\theta,\gamma,V,\nu_0}\{\varrho(v(\cdot), \mu(\cdot)) \leq \delta\} = -I_{\nu,V}(\mu(\cdot)),$$

*where*

$$\tag{3.26}$$
$$I_{\nu,V}(\mu(\cdot)) = \begin{cases} \int_0^T ||\dot{\mu}(s) - R(\mu(s)) - A^*(\mu(s))||_{\mu(s)}^2 \, ds, & \text{if } \mu(0) = \nu \\ & \text{and } \mu(\cdot) \text{ is absolutely continuous,} \\ \infty. & \text{elsewhere} \end{cases}$$

**Remark:** This yields the weak LDP, that is, the upper bound holds for compact sets only.

**Proof:** Theorem 3.4 combined with Theorem 3.5 in [8], we have for $\mu(\cdot) \in C([0,T], M_1^{pro}(E))$,

$$\lim_{\delta \to 0} \liminf_{\gamma \to 0} \gamma \log P_\nu^{\theta,\gamma,\nu_0}\{v : \varrho(v(\cdot), \mu(\cdot)) < \delta\} \tag{3.27}$$
$$= \lim_{\delta \to 0} \limsup_{\gamma \to 0} \gamma \log P_\nu^{\theta,\gamma,\nu_0}\{v : \varrho(v(\cdot), \mu(\cdot)) \leq \delta\} = -I_\nu(\mu(\cdot)).$$

Let

$$C = \sup_{\mu(\cdot) \in C^w([0,T], M_1^{pro}(E))} |\int_0^T \int_E \int_E [\int_E V(x,z)\mu(s,dz)][\int_E V(y,z)\mu(s,dz)]Q(\mu(s);dx,dy)ds|.$$

Since $V$ is uniformly bounded, $C$ is a finite constant. For any measurable subset $B$ of space $C^w([0,T], M_1^{pro}(E))$, by using (3.23), Hölder's inequality, and martingle property, we get for any $\alpha > 0, \beta > 0, \frac{1}{\alpha} + \frac{1}{\beta} = 1$,

$$P_\nu^{\theta,\gamma,V,\nu_0}\{B\} = \int_B Z_V(T) dP_\nu^{\theta,\gamma,\nu_0}$$
$$\leq e^{\frac{C}{2\gamma}} (\int \exp[\frac{\alpha}{\gamma} \int_0^T \int_E (\int_E V(y,z)\mu(s,dz))M(ds,dy)] dP_\nu^{\theta,\gamma,\nu_0})^{1/\alpha} P_\nu^{\theta,\gamma,\nu_0}\{B\}^{1/\beta} \tag{3.28}$$
$$\leq e^{\frac{C}{2\gamma}(1+\alpha)} P_\nu^{\theta,\gamma,\nu_0}\{B\}^{1/\beta} (\int Z_{\alpha V}(T) dP_\nu^{\theta,\gamma,\nu_0})^{1/\alpha}$$
$$= e^{\frac{C}{2\gamma}(1+\alpha)} P_\nu^{\theta,\gamma,\nu_0}\{B\}^{1/\beta}.$$

By choosing $B = \{v : \varrho(v(\cdot), \mu(\cdot)) \leq \delta\}$ in (3.28), we get

$$\lim_{\delta \to 0} \limsup_{\gamma \to 0} \gamma \log P_\nu^{\theta,\gamma,V,\nu_0}\{\varrho(v(\cdot), \mu(\cdot)) \leq \delta\} \leq \frac{C(1+\alpha)}{2} - \frac{1}{\beta} I_\nu(\mu(\cdot)) \tag{3.29}$$

which implies (3.25) for paths that are not absolutely continuous.

Let $\bar{\mu}(\cdot)$ be an arbitrary absolutely continuous path in $C^w([0,T], M_1^{pro}(E))$. For any $\mu \in M_1(E)$, let $V(\mu)(y) = V(\mu, y) = \int_E V(y,z)\mu(dz)$. Then we have

$$\int_0^T \int_E V(\bar{\mu}(s), y) M(ds, dy) = \langle \mu(T), V(\bar{\mu}(T)) \rangle - \langle \mu(0), V(\bar{\mu}(0)) \rangle \tag{3.30}$$
$$- \int_0^T \langle \mu(s), (\frac{\partial}{\partial s} + A) V(\bar{\mu}(s)) \rangle ds,$$

which is continuous in $\mu(\cdot)$ in the topology generated by $\varrho$. On the other hand

$$\int_E \int_E V(\mu, x) V(\mu, y) Q(\mu; dx, dy) = \int_E \int_E \int_E V(x, z_1) V(x, z_2) \mu(dz_1) \mu(dz_2) \mu(dx)$$
$$- \int_E \int_E \int_E \int_E V(x, z_1) V(y, z_2) \mu(dz_1) \mu(dz_2) \mu(dx) \mu(dy),$$

which is also continuous in the topology generated by $\varrho$.

Define

$$\Gamma_V(\mu(\cdot)) = \langle \mu(T), V(\mu(T)) \rangle - \langle \mu(0), V(\mu(0)) \rangle \tag{3.31}$$
$$- \int_0^T \langle \mu(s), (\frac{\partial}{\partial s} + A)V(\mu(s)) \rangle ds$$
$$- \frac{1}{2} \int_0^T \int_E \int_E [\int_E V(x,z)\mu(s,dz)][\int_E V(y,z)\mu(s,dz)]Q(\mu(s); dx, dy)ds.$$

Let $B^\circ(\bar{\mu}(\cdot), \delta)$ denote the interior of $B(\bar{\mu}(\cdot), \delta)$. For any $\varepsilon > 0$, we can choose $\delta > 0$ small enough such that

$$|\liminf_{\gamma \to 0} \gamma \log P_\nu^{\theta, \gamma, \nu_0} \{B^\circ(\bar{\mu}(\cdot), \delta)\} + I_\nu(\bar{\mu}(\cdot))| < \varepsilon, \tag{3.32}$$

$$|\limsup_{\gamma \to 0} \gamma \log P_\nu^{\theta, \gamma, \nu_0} \{B(\bar{\mu}(\cdot), \delta)\} + I_\nu(\bar{\mu}(\cdot))| < \varepsilon, \tag{3.33}$$

$$|\int_0^T \int_E V(\bar{\mu}(s), y) M(ds, dy) - (\langle \bar{\mu}(T), V(\bar{\mu}(T)) \rangle - \langle \bar{\mu}(0), V(\bar{\mu}(0)) \rangle \tag{3.34}$$
$$- \int_0^T \langle \bar{\mu}(s), (\frac{\partial}{\partial s} + A)V(\bar{\mu}(s)) \rangle ds| < \varepsilon,$$

$$|\int_0^T \int_E \int_E V(\mu(s), x) V(\mu(s), y) Q(\mu(s); dx, dy) ds \tag{3.35}$$
$$- \int_0^T \int_E \int_E V(\bar{\mu}(s), x) V(\bar{\mu}(s), y) Q(\mu(s); dx, dy) ds| < \varepsilon,$$

$$|\int_0^T \int_E \int_E V(\bar{\mu}(s), x) V(\bar{\mu}(s), y) Q(\mu; dx, dy) ds \tag{3.36}$$
$$- \int_0^T \int_E \int_E V(\bar{\mu}(s), x) V(\bar{\mu}(s), y) Q(\bar{\mu}(s); dx, dy) ds| < \varepsilon,$$

$$|\int_0^T \int_E \int_E (V(\mu(s), x) - V(\bar{\mu}(s), x))(V(\mu(s), y) - V(\bar{\mu}(s), y)) Q(\mu(s); dx, dy) ds| < \varepsilon. \tag{3.37}$$

By (3.34)-(3.37), we get

$$P_\nu^{\theta, \gamma, V, \nu_0}\{B(\bar{\mu}(\cdot), \delta)\} = \int_{B(\bar{\mu}(\cdot), \delta)} \exp[\frac{1}{\gamma} G_V(\mu(\cdot))] dP_\nu^{\theta, \gamma, \nu_0}$$
$$\leq \exp\{\frac{1}{\gamma}[\Gamma_V(\bar{\mu}(\cdot)) + \frac{3}{2}\varepsilon]\}$$

$$\times \int_{B(\bar{\mu}(\cdot),\delta)} \exp\{\frac{1}{\gamma}[\int_0^T \int_E [\int_E V(y,z)(\mu(ds,dz) - \bar{\mu}(s,dz))]M(ds,dy)]\}dP_\nu^{\theta,\gamma,\nu_0}$$

$$\leq \exp\{\frac{1}{\gamma}[\Gamma_V(\bar{\mu}(\cdot)) + \frac{3+\alpha}{2}\varepsilon]\}P_\nu^{\theta,\gamma,\nu_0}\{B(\bar{\mu}(\cdot),\delta)\}^{1/\beta} \quad (3.38)$$

$$\times (\int \exp\{\frac{1}{\gamma}[\int_0^T \int_E [\int_E \alpha V(y,z)(\mu(ds,dz) - \bar{\mu}(s,dz))]M(ds,dy)$$

$$-\frac{\alpha^2}{2}\int_0^T \int_E \int_E (V(\mu(s),x) - V(\bar{\mu}(s),x))$$

$$\times (V(\mu(s),y) - V(\bar{\mu}(s),y))Q(\mu(s);dx,dy)ds]\}dP_\nu^{\theta,\gamma,\nu_0})^{1/\alpha}$$

$$= \exp\{\frac{1}{\gamma}[\Gamma_V(\bar{\mu}(\cdot)) + \frac{3+\alpha}{2}\varepsilon]\}P_\nu^{\theta,\gamma,\nu_0}\{B(\bar{\mu}(\cdot),\delta)\}^{1/\beta}.$$

Together with (3.33) this implies

$$\lim_{\delta \to 0} \limsup_{\gamma \to 0} \gamma \log P_\nu^{\theta,\gamma,\nu_0}\{B(\bar{\mu}(\cdot),\delta)\} \leq \Gamma_V(\bar{\mu}(\cdot)) + \frac{3+\alpha}{2}\varepsilon - \frac{1}{\beta}I_\nu(\bar{\mu}(\cdot)).$$

By letting $\varepsilon \to 0$, then $\beta \to 1$, we get

$$\lim_{\delta \to 0} \limsup_{\gamma \to 0} \gamma \log P_\nu^{\theta,\gamma,\nu_0}\{B(\bar{\mu}(\cdot),\delta)\} \leq \Gamma_V(\bar{\mu}(\cdot)) - I_\nu(\bar{\mu}(\cdot)). \quad (3.39)$$

By an argument similar to that used in the derivation of (3.38) and Hölder's inequality we get

$$P_\nu^{\theta,\gamma,V,\nu_0}\{B^\circ(\bar{\mu}(\cdot),\delta)\} = \int_{B^\circ(\bar{\mu}(\cdot),\delta)} \exp[\frac{1}{\gamma}G_V(\mu(\cdot))]dP_\nu^{\theta,\gamma,\nu_0}$$

$$\geq \exp\{\frac{1}{\gamma}[\Gamma_V(\bar{\mu}(\cdot)) - \frac{3}{2}\varepsilon]\}$$

$$\times \int_{B^\circ(\bar{\mu}(\cdot),\delta)} \exp\{\frac{1}{\gamma}[\int_0^T \int_E [\int_E V(y,z)(\mu(ds,dz) - \bar{\mu}(s,dz))]M(ds,dy)]\}dP_\nu^{\theta,\gamma,\nu_0}$$

$$\geq \exp\{\frac{1}{\gamma}[\Gamma_V(\bar{\mu}(\cdot)) - (\frac{3}{2} + \frac{\alpha}{2\beta})\varepsilon]\}$$

$$\times \int_{B^\circ(\bar{\mu}(\cdot),\delta)} \exp\{\frac{1}{\gamma}[\int_0^T \int_E (\int_E V(y,z)(\mu(ds,dz) - \bar{\mu}(s,dz)))M(ds,dy)$$

$$+\frac{\alpha}{2\beta}\int_0^T \int_E \int_E (V(\mu(s),x) - V(\bar{\mu}(s),x))$$

$$\times (V(\mu(s),y) - V(\bar{\mu}(s),y))Q(\mu(s);dx,dy)ds]\}dP_\nu^{\theta,\gamma,\nu_0}$$

$$\geq \exp\{\frac{1}{\gamma}[\Gamma_V(\bar{\mu}(\cdot)) - (\frac{3}{2} + \frac{\alpha}{2\beta})\varepsilon]\}P_\nu^{\theta,\gamma,\nu_0}\{B^\circ(\bar{\mu}(\cdot),\delta)\}^\beta \quad (3.40)$$

$$\times (\int \exp\{\frac{1}{\gamma}[\int_0^T \int_E [\int_E -\frac{\alpha}{\beta}V(y,z)(\mu(ds,dz) - \bar{\mu}(s,dz))]M(ds,dy)$$

$$-\frac{\alpha^2}{2\beta^2}\int_0^T \int_E \int_E (V(\mu(s),x) - V(\bar{\mu}(s),x))$$

$$\times (V(\mu(s),y) - V(\bar{\mu}(s),y))Q(\mu(s);dx,dy)ds]\}dP_\nu^{\theta,\gamma,\nu_0})^{-\frac{\beta}{\alpha}}$$

$$= \exp\{\frac{1}{\gamma}[\Gamma_V(\bar{\mu}(\cdot)) - (\frac{3}{2} + \frac{\alpha}{2\beta})\varepsilon]\}P_\nu^{\theta,\gamma,\nu_0}\{B^\circ(\bar{\mu}(\cdot),\delta)\}^\beta,$$

which, combined with (3.32), leads to

$$\lim_{\delta \to 0} \liminf_{\gamma \to 0} \gamma \log P_\nu^{\theta,\gamma,\nu_0}\{B^\circ(\bar{\mu}(\cdot),\delta)\} \geq \Gamma_V(\bar{\mu}(\cdot)) - (\frac{3}{2} + \frac{\alpha}{2\beta})\varepsilon + \beta I_\nu(\bar{\mu}(\cdot)).$$

By letting $\varepsilon \to 0$, then $\beta \to 1$, we get

$$\lim_{\delta \to 0} \liminf_{\gamma \to 0} \gamma \log P_\nu^{\theta,\gamma,\nu_0}\{B^\circ(\bar{\mu}(\cdot),\delta)\} \geq \Gamma_V(\bar{\mu}(\cdot)) - I_\nu(\bar{\mu}(\cdot)). \tag{3.41}$$

The equalities in (3.25) will follow if we can show that for absolutely continuous path $\mu(\cdot)$,

$$I_{\nu,V}(\mu(\cdot)) = -\Gamma_V(\mu(\cdot)) + I_\nu(\mu(\cdot)). \tag{3.42}$$

By direct calculation we have

$$\begin{aligned} I_\nu(\mu(\cdot)) - \Gamma_V(\mu(\cdot)) &= \sup_{f \in C^{1,0}([0,T] \times E)} \{\langle \mu(T), f(T) - V(\mu(T))\rangle - \langle \mu(0), f(0) - V(\mu(0))\rangle \\ &\quad - \int_0^T [\langle \mu(s), (\frac{\partial}{\partial s} + A)f - V(\mu(s))\rangle + \langle R(\mu(s)), f(s) - V(\mu(s))\rangle] \, ds \\ &\quad - \frac{1}{2} \int_0^T \int_E \int_E (f(s,x) - V(\mu(s),x))(f(s,y) - V(\mu(s),y)) \, Q(\mu(s); dx, dy) \, ds\} \\ &= \sup_{f \in C^{1,0}([0,T] \times E)} \{\langle \mu(T), f(T)\rangle - \langle \mu(0), f(0)\rangle \\ &\quad - \int_0^T [\langle \mu(s), (\frac{\partial}{\partial s} + A)f\rangle + \langle R(\mu(s)), f(s)\rangle] \, ds \\ &\quad - \frac{1}{2} \int_0^T \int_E \int_E f(s,x) \, f(s,y) \, Q(\mu(s); dx, dy) \, ds\}, \end{aligned} \tag{3.43}$$

which leads to (3.42).

Finally by (3.42), we get that the level sets associated with $I_{\nu,V}(\mu(\cdot))$ are compact. □

**Acknowledgements.** We thank an anonymous referee for comments and suggestions that helped improve the presentation of this paper.

# Appendix

The proof of Theorem 3.6:

For any $f$ in $C^{1,0}([0,T] \times E)$, and $\mu(\cdot)$ in $\mathcal{H}_\nu$, we have that

$$\langle \mu(T), f(T)\rangle - \langle \mu(0), f(0)\rangle - \int_0^T \langle \mu(u), (\frac{\partial}{\partial u} + A)f\rangle \, du$$

$$-\frac{1}{2}\int_0^T \int_E \int_E f(u,x)\,f(u,y)\,Q(\mu(u);dx,dy)\,du$$

$$= \int_0^T \langle \dot{\mu}(u) - A^*(\mu(u)), f(u)\rangle du - \frac{1}{2}\int_0^T \int_E \int_E f(u,x)\,f(u,y)\,Q(\mu(u);dx,dy)\,du$$

$$\leq \int_0^T ||\dot{\mu}(u) - A^*(\mu(u))||_{\mu(u)}du,$$

which, together with Theorem 3.5, implies that

$$I_\nu(\mu(\cdot)) \leq \int_0^T ||\dot{\mu}(u) - A^*(\mu(u))||_{\mu(u)}du. \tag{A.1}$$

We next show that

$$I_\nu(\mu(\cdot)) \geq \int_0^T ||\dot{\mu}(u) - A^*(\mu(u))||_{\mu(u)}du. \tag{A.2}$$

Since $\nu$ is non-degenerate and $\mu(\cdot)$ is in $\mathcal{H}_\nu$, we have that

$$I_\nu(\mu(\cdot)) = S_\nu(\mu(\cdot)) < \infty.$$

For any $[s,t] \subset [0,T]$, and any $f$ in $C^{1,0}([s,t]\times E)$, let

$$l_{s,t}(f) = \langle \mu(t), f(t)\rangle - \langle \mu(s), f(s)\rangle - \int_s^t \langle \mu(u), (\frac{\partial}{\partial u} + A)f\rangle du.$$

Let

$$L^2([s,t]\times E) = \{f : \int_s^t \int_E f^2(u,x)\mu(u,dx)du < \infty\},$$

and $L$ is the linear subspace of $L^2([s,t]\times E)$ of all functions which are constant in space variable $x$. Let $L^2([s,t]\times E)/L$ be the quotient space of $L^2([s,t]\times E)$ module $L$. We introduce on $L^2([s,t]\times E)/L$ the following norm

$$||h|| = \Big\{\int_s^t \int_E \int_E h(u,x)h(u,y)Q(\mu(u);dx,dy)du\Big\}^{1/2},$$

and the inner product

$$\langle h_1, h_2\rangle = \int_s^t \int_E \int_E h_1(u,x)h_2(u,y)Q(\mu(u);dx,dy)du.$$

It is not hard to check that space $(L^2([s,t]\times E)/L, ||\cdot||)$ is a pre-Hilbert space. Let $\mathcal{L}^2([s,t]\times E)$ be the completion of space $L^2([s,t]\times E)/L$ (cf. for existence refer to page 56 of Yosida [15]). Then $(\mathcal{L}^2([s,t]\times E), ||\cdot||)$ becomes a Hilbert space.

Let $\mathcal{L}^2_{sub}([s,t]\times E)$ denote the closure in $\mathcal{L}^2([s,t]\times E)$ of the linear span of the set $C^{1,0}([s,t]\times E)$. By an argument similar to Dawson and Gärtner [3] on pages 277-280, and the Hahn-Banach Extension Theorem (cf. page 106 of Yosida [15]) we have that $l_{s,t}$ is a bounded linear functional on $\mathcal{L}^2_{sub}([s,t] \times E)$. Thus by the Riesz Representation Theorem, one can find a function $h$ in $\mathcal{L}^2_{sub}([0,T] \times E)$ such that

$$l_{s,t}(f) = \int_s^t \int_E \int_E f(u,x)h(u,y)Q(\mu(u);dx,dy)du, \tag{A.3}$$

$$\inf_{f \in C^{1,0}([s,t] \times E)} \int_s^t \int_E \int_E (h(u,x) - f(u,x))(h(u,y) - f(u,y))Q(\mu(u);dx,dy)du = 0, \tag{A.4}$$

$$S_\nu(\mu(\cdot)) = \frac{1}{2} \int_0^T \int_E \int_E h(u,x)h(u,y)Q(\mu(u);dx,dy)du, \tag{A.5}$$

which implies that

$$\int_0^T ||\dot{\mu}(u) - A^*(\mu(u))||_{\mu(u)} du = \int_0^T \sup_{g \in C^\infty(E)} \{\langle \dot{\mu}(u) - A^*, g \rangle$$
$$- \frac{1}{2} \int_E \int_E g(x)g(y)Q(\mu(u);dx,dy)\}$$
$$= \int_0^T \sup_{f \in C^{1,0}([0,T] \times E)} \{\langle \dot{\mu}(u) - A^*, f(u) \rangle - \langle \mu(u), \dot{f}(u) \rangle$$
$$- \frac{1}{2} \int_E \int_E f(u,x)f(u,y)Q(\mu(u);dx,dy)\}$$
$$= S_\nu(\mu(\cdot)) - \frac{1}{2} \inf_{f \in C^{1,0}([0,T] \times E)} \int_0^T \int_E \int_E (h(u,x) - f(u,x))(h(u,y) - f(u,y))Q(\mu(u);dx,dy)du$$
$$= S_\nu(\mu(\cdot)),$$

where $(A.3)$ is used in deriving the third equality, and $(A.4)$ and $(A.5)$ are used in reaching the last equality. This together with Theorem 3.5 implies (3.22). □

# References

[1] Bauer,H. *Probability Theory and Elements of Measure Theory*. Academic Press, London 1981.

[2] Dawson, D.A.(1978). Geostochastic calculus. *Canadian Journal of Statistics*, **6**:143–168.

[3] Dawson, D.A. and Gärtner, J.(1987). Large deviations from the McKean-Vlasov limit for weakly interacting diffusions. *Stochastics*, **20**:247–308.

[4] Dembo, A. and Zeitouni, O. *Large Deviations and Applications*. Jones and Bartlett Publishers, Boston 1993.

[5] Donsker, M.D. and Varadhan, S.R.S.(1975). Asymptotic evaluation of certain Markov process expectations for large time, I. *Comm. Pure Appl. Math.*, **28**, pp. 1-47.

[6] Ethier, S.N. and Kurtz, T.G.(1993). Fleming-Viot processes in population genetics. *SIAM J. Control and Optimization*, **31**, N.2:345–386.

[7] Ethier, S.N. and Kurtz, T.G.(1994). Convergence to Fleming-Viot processes in the weak atomic topology. *Stoch. Proc. Appl.*, **54**, 1–27.

[8] Freidlin, M.I. and Wentzell, A.D. *Random Perturbations of Dynamical Systems.* Springer-Verlag, New York 1984.

[9] Hartl, D.L. and Clark, A.G. *Principles of Population Genetics*, 2nd ed. Sinauer Associates, Inc., Sunderland, Massachusetts, 1989.

[10] Pukhalskii, A.A.(1991). On functional principle of large deviations. in New Trends in Probability and Statistics, ed. V. Sazonov and T. Shervashidze, VSP Moks'las, Moskva.

[11] Shiga, T.(1981). Diffusion processes in population genetics. *J. Math. and Kyoto Univ.*, **21**:133–151.

[12] Shiga, T.(1990). A stochastic equation based on a Poisson system for a class of measure-valued diffusion processes. *J. Math. and Kyoto Univ.*, **30**:245–279.

[13] Stroock, D.W. and Varadhan, S.R.S. *Multidimensional Diffusion Processes.* Springer-Verlag, Berlin 1979.

[14] Wright, S.(1949). Adaptation and selection. in *Genetics, Paleontology, and Evolution.* ed. G.L. Jepson et al, 365-389, Princeton Univ. Press.

[15] Yosida, K. *Functional Analysis, Sixth Ed..* Springer-Verlag, Berlin 1980.

Donald A. Dawson
The Fields Institute
222 College Street
Toronto, Ontario
Canada M5T 3J1
e-mail: don@fields.fields.utoronto.ca

Shui Feng
Department of Mathematics
and Statistics
McMaster University
Hamilton, Ontario
Canada L8S 4K1
e-mail:shuifeng@mcmail.cis.mcmaster.ca